\documentclass[12pt]{article}
\usepackage[latin1]{inputenc}
\usepackage{amsmath}
\usepackage{enumerate}
\usepackage{amsfonts}
\usepackage{amssymb}

\newcommand{\N}{\mathbb{N}}

\newcommand{\RP}{\mathbb{R}^{{\scriptscriptstyle{+}}}}

\newcommand{\E}{\mathcal{E}}

\newcommand{\B}{\mathcal{B}}

\newcommand{\A}{\mathcal{A}}

\newcommand{\C}{\mathcal{C}}

\newcommand{\D}{\mathcal{D}}
\newcommand{\U}{\mathcal{U}}
\newcommand{\F}{\mathcal{F}}

\newcommand{\Av}{\mathrm{Av}}

\newcommand{\card}{\mathrm{card}}

\newcommand{\vtau}{\vec{\tau}}

\newcommand{\Vartriangle}{\vartriangle}

\newcommand{\lo}{\mathcal{L}_0}

\newcommand{\vz}{\vec{z}}

\newcommand{\vy}{\vec{y}}

\newcommand{\vu}{\vec{u}}

\newcommand{\imin}[1]{#1^{-1}}

\newcommand{\qed}{\hfill $\Box$ \vspace{\baselineskip}}

\renewcommand{\dim}{\textbf{Proof.}}

\numberwithin{equation}{section}

\newtheorem{teo}{Theorem}[section]

\newtheorem{df}[teo]{Definition}

\newtheorem{cor}[teo]{Corollary}

\newtheorem{prop}[teo]{Proposition}

\newtheorem{lem}[teo]{Lemma}

\newcommand{\be}{\begin{equation}}
\newcommand{\ee}{\end{equation}}

\begin{document}

 \noindent{\bf Title: }         Maharam's problem\\
 {\bf  Author: }        Michel Talagrand\\
{\bf Subj-class:}       FA\\
 {\bf Comments:}       30p.\\
 {\bf MSC-class:}      28A12.
\bigskip\begin{center}
{\huge Maharam's problem }
\end{center}

\medskip

\begin{center}
 Michel Talagrand \footnote {Work partially supported by an NSF grant}\ \\
Universit\'e Paris VI and 
the Ohio State University \\
Dedicated to J. W. Roberts
\end{center}

\vspace{1cm}

\begin{abstract} We  construct an exhaustive submeasure that is
not equivalent to a measure. This solves problems of J. von
Neumann (1937) and D. Maharam (1947).
\end{abstract}

\section{Introduction}

Consider a Boolean algebra $\B$ of sets. A map $\nu:\B \to \RP$ is
called a \emph{submeasure} if it satisfies the following
properties:
\begin{eqnarray}
&&   \nu(\varnothing)=0, \\
&&   A \subset B \ \ A,B \in \B \Longrightarrow \nu(A) \leq \nu(B),  \\
&&   A,B \in \B \Rightarrow \nu(A \cup B) \leq \nu(A) + \nu(B).
\end{eqnarray}
If we have $\nu(A \cup B)=\nu(A)+\nu(B)$ whenever $A$ and $B$ are
disjoint, we say that $\nu$ is a (finitely additive)
\emph{measure}.

We say that a sequence $(E_n)$ of $\B$ is \emph{disjoint} if $E_n
\cap E_m= \varnothing$ whenever $n \neq m$. A submeasure
  is \emph{exhaustive} if $\lim_{n
\to \infty}\nu(E_n)=0$ whenever $(E_n)$ is a disjoint sequence of
$\B$. A measure is obviously exhaustive. Given two submeasures
$\nu_1$ and $\nu_2$, we say that $\nu_1$ is \emph{absolutely
continuous with respect to} $\nu_2$ if
\begin{equation}\label{1.4}
\forall \varepsilon >0,\, \exists \alpha >0, \, \nu_2(A)\leq \alpha
\Longrightarrow\nu_1(A) \leq \varepsilon.
\end{equation}

If a submeasure is absolutely continuous with respect to a
measure, it is exhaustive. One of the many equivalent forms of
Maharam's problem is whether the converse is true.\\
\emph{Maharam's problem:} If a submeasure is exhaustive, is it
absolutely continuous with respect to a measure?

In words, we are asking whether the only way a submeasure can be
exhaustive is because it really resembles a measure. This question
has been one of the longest standing classical questions of
measure theory. It occurs in a variety of forms (some of which
will be discussed below).

Several important contributions were made to Maharam's problem. N.
Kalton and J. W. Roberts proved \cite{K-R} that a submeasure is
absolutely continuous with respect to a measure if (and, of
course, only if) it is uniformly exhaustive, i.e.
\begin{equation}\label{1.5}
\forall \varepsilon >0, \ \ \exists n, \ \ E_1,\ldots,E_n \,\,
\text{disjoint} \Longrightarrow \inf_{i \leq n} \nu(E_i) \leq
\varepsilon.
\end{equation}
Thus Maharam's problem can be reformulated as to whether an
exhaustive submeasure is necessarily uniformly exhaustive. Two other
fundamental contributions by J.W. Roberts \cite{R} and I. Farah
\cite{F} are used in an essential way in this paper and will be
discussed in great detail later.

We prove that Maharam's problem has negative answer.

\begin{teo}\label{teo1.1} There exists a non-zero exhaustive submeasure $\nu$ on
the algebra $\B$ of clopen sets of the Cantor set that is not
uniformly exhaustive (and thus is not absolutely continuous with
respect to a measure). Moreover, no non-zero measure $\mu$ on $\B$
is absolutely continuous with respect to $\nu$.
\end{teo}

We now spell out some consequences of Theorem \ref{teo1.1}. It has
been known for a while how to deduce these results from Theorem
\ref{teo1.1}. For the convenience of the reader these (easy)
arguments will be given in a self-contained way in the last
section of the paper.

Since Maharam's original question and von Neumann problem are
formulated in terms of general Boolean algebras (i.e., that are
not a priori represented  as algebras of sets) we must briefly
mention these. We will denote by $0$ and $1$ respectively the
smallest and the largest element of a Boolean algebra $\B$, but we will
denote the operations by $\cap,\cup$, etc. as in the case of
algebras of sets. A Boolean algebra $\B$ is called $\sigma$-complete if
any countable set $\C$ has a least upper bound $\cup \C$ (and thus a
greatest lower bound $\cap \C$). A submeasure $\nu$ on $\B$ is
called \emph{continuous} if whenever $(A_n)$ is a decreasing
sequence with $\bigcap_nA_n=0$ we have $\lim_{n \to \infty}
\nu(A_n)=0$. The submeasure is called \emph{positive} if $\nu(A)=0
\Longrightarrow A=0$.

A $\sigma$-complete algebra $\B$ on which there is a positive
continuous submeasure is called a \emph{submeasure algebra}. If
there is a positive continuous measure on $\B$, $\B$ is called a
\emph{measure algebra}.

Probably the most important consequence of our construction is that it proves the existence of radically new Boolean algebras

\begin{teo}\label{teo1.2} There exists a submeasure algebra $\B$
that is not a measure algebra. In fact, there does not exist a
positive measure on $\B$, and there not exist  a (non-zero) continuous
measure on $\B$ either.
\end{teo}

This answers a question raised by D. Maharam in her 1947 paper
\cite{M}.

A subset $\C$ of a boolean algebra $\B$ is called \emph{disjoint}
if $A \cap B =\varnothing$ whenever $A,B \in \C$, $A \neq B$. A
disjoint set $\C$ is called a \emph{partition} if $\cup \C=1$ (= the
largest element of $\B$). If every disjoint collection of $\B$ is
countable, $\B$ is said to satisfy the \emph{countable chain condition}.

If $\Pi$ is a partition of $\B$ we say that $A \in \B$ is
\emph{finitely covered} by $\Pi$ if there is a finite subset
$\{A_1,\ldots,A_n\}$ of $\Pi$ with $A \subset \bigcup_{i \leq n}
A_i$. We say that $\B$ satisfies the \emph{general distributive
law} if whenever $(\Pi_n)$ is a sequence of partitions of $\B$,
there is a single partition $\Pi$ of $\B$ such that every element of
$\Pi$ is finitely covered by each $\Pi_n$. (This terminology is not used by every author, such a $\sigma$-algebra is called weakly $(\sigma-\infty)$ distributive in \cite{F1}.)

\begin{teo}\label{teo1.3} There exists a $\sigma$-complete algebra
that satisfies the countable chain condition and the general
distributive law, but is not a measure algebra.
\end{teo}

We spell out this statement because it answers negatively a
problem raised by J. von Neumann in the Scottish book (\cite{Ma}
problem 163), but it is a simple consequence of Theorem
\ref{teo1.2}, since every submeasure algebra satisfies the
countable chain condition and the general distributive law.
Examples of this type had been known under special axioms, such as
the negation of Suslin's hypothesis \cite{M}, but our example is the first one that does not use any special axiom. (In fact, it has been recently shown \cite{bjp}, \cite{V}, that essentially the only way to produce a counterexample to von Neumann problem that does not use special axioms is indeed to solve Maharam's problem.) 

Consider now a topological vector space $X$ with a metrizable
topology, and $d$ a translation invariant distance that defines
this topology. If $\B$ is a boolean algebra of subsets of a set
$T$, an ($X$-valued) \emph{vector measure} is a map $\theta:\B \to
X$ such that $\theta(A \cup B)=\theta(A)+\theta(B)$ whenever $A
\cap B=\varnothing$. We say that it is \emph{exhaustive} if  $\lim_{n \to
\infty} \theta(E_n)=0$ for each disjoint sequence $(E_n)$ of $\B$.
A positive measure $\mu$ on $\B$ is called a \emph{control
measure} for $\theta$ if
$$\forall \varepsilon >0,\ \ \exists \alpha >0, \mu(A) \leq \alpha
\Longrightarrow d(0,\theta(A))\leq \varepsilon.$$

\begin{teo}\label{teo1.4}{\bf  (Negative solution to the Control Measure Problem)}
There exists an exhaustive vector-valued measure that does not
have a control measure.
\end{teo}

We now explain the organization of the paper. The submeasure we
will construct is an object of a rather new nature, since it is
very far from being a measure. It is unlikely that a very simple
example exists at all, and it should not come as a surprise that
our construction is somewhat involved. Therefore it seems
necessary to explain first the main ingredients on which the
construction relies. The fundamental idea is due to J. W. Roberts
\cite{R} and is detailed in Section \ref{2}. Another crucial part
of the construction is a technical device invented by I. Farah
\cite{F}. In Section \ref{3}, we produce a kind of ``miniature
version'' of Theorem \ref{teo1.1}, to explain Farah's device, as
well as some of the other main ideas. The construction of $\nu$
itself is given in Section \ref{4}, and the technical work of
proving that $\nu$ is not zero and is exhaustive is done in
Sections \ref{5} and \ref{6} respectively. Finally, in Section
\ref{7} we give the simple (and known) arguments needed to deduce
Theorems \ref{teo1.2} to \ref{teo1.4} from Theorem \ref{teo1.1}.\\

\noindent {\bf Acknowledgments.} My warmest thanks go to I. Farah
who explained to me the importance of Roberts's work \cite{R},
provided a copy of this hard-to-find paper, rekindled my interest
in this problem, and, above all, made an essential technical
contribution without which my own efforts could hardly have
succeeded.

\section{Roberts}\label{2}

Throughout the paper we write
$$T=\prod_{n \geq 1} \{1,\ldots,2^n\}.$$

\noindent For $\vz \in T$, we thus have $\vz=(z_n)$, $z_n \in
\{1,\ldots,2^n\}$. We denote by $\B_n$ the algebra generated by
the coordinates of rank $\leq n$, and $\B=\bigcup_{n \geq 1}\B_n$
the algebra of the clopen sets of $T$. It is isomorphic to the
algebra of the clopen sets of the Cantor set $\{0,1\}^\N$.

We denote by $\A_n$ the set of atoms of $\B_n$. These are sets of
the form
\begin{equation}\label{2.0}
\{\vz \in T;\, z_1=\tau_1,\ldots,z_n=\tau_n\}
\end{equation}
where $\tau_i$ is an integer $\leq 2^i$. An element $A$ of $\A_n$
will be called an atom of rank $n$.

\begin{df}\cite{R} Consider $1 \leq m<n$. We say that a subset $X$
of $T$ is $(m,n)$-thin if
$$\forall A \in \A_m,\ \ \exists A' \in \A_n,\ \ A' \subset A,\ \
A' \cap X=\varnothing.$$
\end{df}

In words, in each atom of rank $m$, $X$ has a hole big enough to
contain an atom of rank $n$. It is obvious that if $X$ is
$(m,n)$-thin, it is also $(m,n')$-thin when $n'\geq n$.

\begin{df}\cite{R} Consider a (finite) subset $I$ of $\N^*=\N \setminus \{0\}$. We say
that $X \subset T$ is $I$-thin if $X$ is $(m,n)$-thin whenever
$m<n$, $m,n \in I$.
\end{df}

 We denote by $\card I$ the cardinality of a
finite set $I$. For two finite sets $I,J \subset \N^*$, we write
$I\prec J$ if $\max I \leq \min J$.

The following is implicit in \cite{R} and explicit in \cite{F}.

\begin{lem}\label{lem2.3}{\bf (Roberts's selection lemma).} Consider
two integers $s$ and $t$, and sets $I_1,\ldots,I_s \subset \N^*$ with
$\card I_\ell \geq st$ for $1\leq \ell \leq s$. Then we can
relabel the sets $I_1,\ldots,I_s$ so that we can find sets $J_\ell
\subset I_\ell$ with $\card J_\ell=t$ and $J_1\prec J_2\prec
\cdots \prec J_s$.
\end{lem}
\dim\ \ Let us enumerate $I_\ell=\{i_{1,\ell},\ldots,i_{st,\ell}\}$
where $i_{a,\ell}<i_{b,\ell}$ if $a<b$. We can relabel the sets
$I_\ell$ in order to ensure that
\begin{eqnarray*}
\forall k \geq 1 && i_{t,1} \leq i_{t,k} \\
\forall k \geq 2 && i_{2t,2} \leq i_{2t,k}
\end{eqnarray*}
and more generally, for any $\ell<s$ that
\begin{eqnarray}\label{2.1}
\forall k \geq \ell && i_{\ell t,\ell} \leq i_{\ell t,k}
\end{eqnarray}
We then define
$$J_\ell=\{i_{(\ell-1)t+1,\ell},\ldots,i_{\ell t,\ell}\}.$$
To see that for $1 \leq \ell <s$ we have $J_\ell \prec J_{\ell+1}$
we use \eqref{2.1} for $k=\ell+1$, so that $i_{\ell t,\ell} \leq
i_{\ell t,\ell+1}<i_{\ell t+1,\ell+1}$. \qed

The reader might observe that it would in fact suffice to assume
that $\card I_\ell \geq s(t-1)+1$; but this refinement yields no
benefits for our purposes.

Throughout the paper, given an integer $\tau \leq 2^n$, we write
\begin{equation}\label{2.2}
S_{n,\tau}=\{\vz \in T;\,z_n \neq \tau\}
\end{equation}
so that its complement $S^c_{n,\tau}$ is the set $\{ \vz\in T;\,
z_n=\tau\}$. Thus on the set $S_{n,\tau}$ we forbid that the ${\rm n}^{{\rm th}}$ coordinate of $\vz$ be $\tau$ while on $S^c_{n,\tau}$ we force it to be $\tau$. 

\begin{prop}\label{prop2.4} Consider sets $X_1,\ldots,X_q \subset T$, and assume
that for each $\ell \leq q$ the set $X_\ell$ is $I_\ell$-thin, for
a certain set $I_\ell$ with $\card I_\ell \geq 3q$. Then for each
$n$ and each integer $\tau \leq 2^n$ we have
\begin{equation}\label{2.3}
S_{n,\tau}^c \not\subset \bigcup_{\ell \leq q}X_\ell.
\end{equation}
\end{prop}
\dim\ \ We use Lemma \ref{lem2.3} for $s=q$ and $t=3$ to produce
sets $J_\ell \subset I_\ell$ with $J_1\prec J_2\prec \cdots \prec
J_q$ and $\card J_\ell = 3$. Let $J_\ell=(m_\ell,n_\ell,r_\ell)$, so $r_\ell \leq
m_{\ell+1}$ since $J_\ell\prec J_{\ell+1}$.

To explain the idea (on which the paper ultimately relies) let us
prove first that $T \not \subset \bigcup_{\ell \leq q}X_\ell$. We make
an inductive construction to avoid in turn the sets $X_\ell$. We
start with any $A_1 \in \A_{m_1}$. Since $X_1$ is
$(m_1,n_1)$-thin, we can find $C_1 \in \A_{n_1}$ with $ C_1 \subset A_1$ and $C_1 \cap
X_1 = \varnothing$. Since $n_1 \leq m_2$ we can find $A_2 \in
\A_{m_2}$ and $A_2 \subset C_1$, and we continue in this manner.
The set $C_q$ does not meet any of the sets $X_\ell$.

To prove  \eqref{2.3}, we must ensure that $C_q \cap
S_{n,\tau}^c \neq \varnothing$. The fundamental fact is that at
each stage we have \emph{two} chances to avoid $X_\ell$, using
either that $X_\ell$ is $(m_\ell,n_\ell)$-thin or that it is
$(n_\ell,r_\ell)$-thin. The details of the construction depend on
the ``position'' of $n$ with respect to the sets $J_\ell$. Rather
that enumerating the cases, we explain what happens when $m_1 <n
\leq r_1$, and this should make what to do in the other cases
obvious.\\
\emph{Case 1.} We have $m_1<n\leq n_1$. Since $S_{n,\tau}\in
\B_n \subset \B_{n_1}$, we can choose $A_1 \in \A_{n_1}$ with $A_1
\subset S^c_{n,\tau}$. Since $X_1$ is $(n_1,r_1)$-thin, we choose
$C_1 \in \A_{r_1}$ with $C_1 \subset A_1$ and $C_1 \cap X_1 =
\varnothing$. We then continue as before, choosing $A_2 \subset
C_1$, $A_2 \in \A_{m_2}$, etc.\\
\emph{Case 2.} We have $m_1<n_1<n \leq r_1$. We choose any $A_1 \in
\A_{m_1}$. Since $X$ is $(m_1,n_1)$-thin, we can choose $C_1 \in
\A_{n_1}$ with $C_1 \subset A_1$ and $C_1 \cap X_1 =\varnothing$.
It is obvious from \eqref{2.0} that, since $n_1<n$, we have $C_1
\cap S^c_{n,\tau} \neq \varnothing$. Since $C_1 \cap S_{n,\tau}^c
\in \B_n \subset \B_{r_1} \subset \B_{m_2}$, we can find $A_2
\subset C_1 \cap S_{n,\tau}^c$, $A_2 \in \A_{m_2}$, and we
continue as before.
\qed

\begin{df}\label{df2.4} Given $\varepsilon >0$, a submeasure $\nu$
on an algebra $\B$ is called $\varepsilon$-\emph{exhaustive} if for each
disjoint sequence $(E_n)$ of $\B$ we have $\limsup_{n \to
\infty}\nu(E_n) \leq \varepsilon$.
\end{df}

\begin{teo} {\bf (Roberts)} For each $q$ there exists a submeasure
$\nu$ on $T$ such that
\begin{eqnarray}
\forall n,\, \forall \tau \leq 2^n,\ \ \nu(S^c_{n,\tau})=1 \label{2.4}\\
\text{$\nu$ is $\frac{1}{q+1}$-exhaustive.} \label{2.5}
\end{eqnarray}
\end{teo}

Of course, \eqref{2.4} implies that $\nu$ is not uniformly
exhaustive. Let us consider the class $\C$ of subsets $X$ of $T$
that are $I$-thin (for a set $I$ depending on $X$) with $\card I
\geq 3q$. For $B \in \B$ we define
\begin{equation}\label{2.5a}
\nu(B)=\min \left(1,\inf\left\{\frac{1}{q+1} \card F;\, F \subset
\C;\, B \subset \cup F\right\}\right),
\end{equation}
where $F$ runs over the finite subsets of $\C$ and $\cup F$
denotes the union of $F$. It is obvious that $\nu$ is a
submeasure, and \eqref{2.4} is an immediate consequence of
Proposition \ref{prop2.4}.

To prove \eqref{2.5} it suffices, given a disjoint sequence
$(E_n)$ of $\B$, to prove that $\liminf_{n \to \infty}\nu(E_n) \leq
1/(1+q)$.

For $X \subset T$, let us write
\begin{equation}\label{2.5b}
(X)_m= \bigcap \{B \in \B_m;\,B \supset X\}= \bigcup \{A,\, A \in
\A_m,\, A \cap X \neq \varnothing \}.
\end{equation}
Since each algebra $\B_m$ is finite, by taking a subsequence we
can assume that for some integers $m(n)$ we have $E_n \in
\B_{m(n)}$, while
\begin{equation}\label{2.6}
\forall k > n,\ \ (E_k)_{m(n)}=(E_{n+1})_{m(n)}.
\end{equation}

We claim that for each $k>n+1$, $E_{k}$ is $(m(n),m(n+1))$-thin. To
prove this, consider $A \in \A_{m(n)}$. If $A \cap E_{k}=
\varnothing$, any $A' \in \A_{m(n+1)}$ with
$A' \subset A$ satisfies $A' \cap E_{k}= \varnothing$. Otherwise $A
\subset (E_{k})_{m(n)}=(E_{n+1})_{m(n)}$ by \eqref{2.6}.
Therefore, $E_{n+1} \cap A \neq \varnothing$. Since $E_{n+1} \in
\B_{m(n+1)}$, we can find $A' \in \A_{m(n+1)}$ with $A' \subset A$
and $A' \subset E_{n+1}$. But then $A' \cap E_{k} = \varnothing$
since $E_{n+1}$ and $E_{k}$ are disjoint. This proves the claim.

 It follows that  for $n \geq 3q+1$, $E_n$ is $I$-thin for
$I=(n(1),\ldots,n(3q))$ and thus $E_n \in \C$, so that $\nu(E_n)
\leq 1/(q+1).$
\qed

\section{Farah}\label{3}

In \cite{F} I. Farah constructs for each $\varepsilon$ an
$\varepsilon$-exhaustive submeasure $\nu$ that is also
pathological, in the sense that every measure that is absolutely
continuous with respect to $\nu$ is zero. In this paper, we
learned several crucial technical ideas, that are essential for
our approach. The concepts and the techniques required to prove
Proposition \ref{3.5} below are essentially all Farah's.

A class $\C$ of \emph{weighted sets} is a subset of $\B \times
\RP$. For a finite subset $F = \{(X_1,w_1),\ldots, (X_n,w_n)\}$ of $\C$, we write throughout the
paper
\begin{equation}\label{3.1}
w(F)=\sum_{i \leq n}w_i;\ \ \cup F=\bigcup_{i \leq n}X_i,
\end{equation}
and for $B \in \B$ we set
\begin{equation}\label{3.2}
\varphi_\C(B)=\inf \{w(F);\, B \subset \cup F\}.
\end{equation}
It is immediate to check that $\varphi_\C$ is a submeasure
provided $\varphi_\C(T)<\infty$. This construction generalizes
\eqref{2.5a}. It is generic; for a submeasure $\nu$, we have
$\nu=\varphi_\C$ where $\C=\{(B,\nu(B));\, B \in \B\}$. Indeed, it
is obvious that $\varphi_\C \leq \nu$, and the reverse inequality
follows by subadditivity of $\nu$.

For technical reasons, when dealing with classes of weighted sets,
we find it convenient to keep track for each pair $(X,w)$
of a distinguished finite subset $I$ of $\N^*$. For this reason we define a
class of \emph{marked weighted sets} as a subset of $\B \times
\cal F \times \RP$, where $\cal F$ denotes the collection of finite subsets of $\N^*$. 

For typographical convenience we write
\begin{equation}\label{3.3}
\alpha(k)={1 \over (k+5)^3}
\end{equation}
and we consider a sequence $(N(k))$ that be specified later. The
specific choice is anyway completely irrelevant, what matters is
that this sequence increases fast enough. In fact, there is
nothing magic about the choice of $\alpha(k)$ either. Any sequence
such that $\sum_kk\alpha(k)<\infty$ would do. We like to stress
than none of the numerical quantities occurring in our
construction plays an essential role. These are all simple choices
that are made for convenience. No attempts whatsoever  have been
made to make optimal or near optimal choices. Let us also point
out that for the purpose of the present section it would work just
fine to take $\alpha(k)=\imin {(k+5)}$, and that the reasons for
taking a smaller value will become clear only in the next section.
For $k \geq 1$ we define the class $\D_k$ of marked weighted sets
by
\begin{equation}\label{3.4}\begin{split}
\D_k=\left\{(X,I,w);\, \exists (\tau(n))_{n \in I},\, X=\bigcap_{n
\in I} S_{n,\tau(n)};\, \card I \leq N(k),\,\right. \\ \left. w
= 2^{-k}\left({N(k) \over \card I}\right)^{\alpha(k)} \right\}.
\end{split}\end{equation}

 The most important part of $\D_k$ consists of the triplets $(X,I,w)$ where $\card I = N(k) $ and $w= 2^{-k}$. The purpose of the relation 
$w= 2^{-k}\left({N(k) / \card I}\right)^{\alpha(k)}$ is to allow the crucial Lemma \ref{lem3.1} below. To understand the relation between the different classes $\D_k$ it might help to observe the following. Whenever $X$ and $I$ are as in \eqref{3.4} and whenever $N(k) \geq \card I$ we have $(X,I, w_k) \in \D_k$ for $w_k= 2^{-k}\left({N(k) / \card I}\right)^{\alpha(k)}$. If we assume, as we may, that the sequence $ 2^{-k}N(k)^{\alpha(k)}$ increases, we see that the sequence $(w_k)$ increases. It is then the smallest value possible of $k$ that gives the smallest possible value of $w_k$, which, as will be apparent in the formula \eqref{3.6} below is the only one that matters.

Given a subset $J$ of $\N^*$ we say that a subset $X$ of $T$ \emph{depends only on the coordinates of
rank} $n \in J$ if whenever $\vz,\vz' \in T$ are such that $z_n=z'_n$
for $n \in J$, we have $\vz \in T$ iff $\vz ' \in T$.
Equivalently, we sometimes say that such a set does not depend on
the coordinates of rank $n \in J^c=\N^* \setminus J$. One of the key 
ideas of the definition of $\D_k$ is the following simple fact.

\begin{lem}\label{lem3.1} Consider $(X,I,w) \in \D_k$ and $J
\subset \N^*$. Then there is $(X',I',w') \in \D_k$ such that $X
\subset X'$, $X'$ depends only on the coordinates in $J$ and
\begin{equation}\label{3.5}
w'=w\left({\card I \over \card I \cap J}\right)^{\alpha(k)}.
\end{equation}
\end{lem}

Since $\alpha(k)$ is small, $w'$ is not really larger than $w$
unless $\card I \cap J \ll \card I$. In particular, since
$\alpha(k) \leq 1/2$ we have
\begin{equation}\label{3.5a}
\card I \cap J \geq \frac{1}{4} \card I \Longrightarrow w' \leq 2w.
\end{equation}
\dim\ \ We define $(X',I',w')$ by \eqref{3.5} and $I'=I \cap J$,
$$
X'= \bigcap_{n \in I'} S_{n,\tau(n)},
$$
where $\tau(n)$ is as in \eqref{3.4}.
\qed

A class of marked weighted sets is a subset of $\B \times
\cal F \times \RP$. By projection onto $\B \times \RP$, to
each class $\C$ of marked weighted sets, we can associate a class
$\C^*$ of weighted sets. For a class $\C$ of marked weighted sets,
we then define $\varphi_\C$ as $\varphi_{\C^*}$ using \eqref{3.2}.
As there is no risk of confusion, we will not distinguish between
$\C$ and $\C^*$ at the level of notation. We define 
\begin{equation}\label{3.6}
\D=\bigcup_{k \geq 1} \D_k;\ \ \psi=\varphi_\D.
\end{equation}

\begin{prop}\label{prop3.2} Let us assume that
\begin{equation}\label{3.7}
N(k) \geq 2^{k+6}(2^{k+5})^{1/\alpha(k)}.
\end{equation}
Then $\psi(T) \geq 2^5$. Moreover $\psi$ is pathological in the
sense that if a measure $\mu$ on $\B$ is absolutely continuous
with respect to $\psi$, then $\mu=0$.
\end{prop}

 Pathological submeasures seem to have been constructed first implicitly in  \cite{D-R} and explicitly in \cite{P}. 

\dim\ \ To prove that $\psi(T) \geq 2^5$, we consider a finite
subset $F$ of $\D$, with $w(F) < 2^5$, and we prove that $T
\not\subset \cup F$. For $k \geq 1$ let $F_k=F \cap \D_k$. For $(X,I,w) \in
\D_k$, we have $w \geq 2^{-k}$, so that $\card F_k \leq 2^{k+5}$
since $w(F_k) \leq w(F) <2^5$. Also we have
$$
2^{-k}\left({N(k) \over \card I}\right)^{\alpha(k)} = w \leq
w(F) \leq 2^5,
$$
so that $\card I \geq (2^{k+5})^{-1/\alpha(k)}N(k):=c(k)$. Thus, under \eqref{3.7} we have $c(k) \geq 2^{k+6}$. Let us
enumerate $F$ as a sequence $(X_r,I_r,w_r)_{r \leq r_0}$
(where $r_0=\card F$) in such a way that if $(X_r,I_r,w_r) \in F_{k(r)}$,
the sequence $k(r)$ is non-decreasing. Since
$$
\sum_{\ell < k} \card F_\ell \leq \sum_{\ell < k} 2^{\ell+5} <
2^{k+5},
$$
we see that $r \geq 2^{k+5}$ implies $k(r) \geq k$ and thus $\card
I_r \geq c(k)$. If we assume \eqref{3.7} we then see that $\card
I_r \geq r+1$. Indeed this is true if $r<2^6$ because $\card I \geq c(1) \geq 2^6$, and  if $r \geq 2^6$ and if $k$ is the
largest integer with $r \geq 2^{k+5}$, then $c(k) \geq 2^{k+6} \geq r+1$. We can then  pick
inductively integers $i_r \in I_r$ that are all different. If
$X_r=\bigcap_{n \in I_r} S_{n,\tau_r(n)}$, any $\vz$ in $T$ with
$z_{i_r}=\tau_r(i_r)$ for $r \leq r_0$ does not belong to any of
the sets $X_r$, and thus $\cup F \neq T$. This proves  that $\psi(T) \geq 2^5$.

We prove now that $\psi$ is pathological. Consider a measure $\mu$ on $\B$, and assume that there exists $k$ such that
$$
\psi(B) \leq 2^{-k} \Longrightarrow \mu(B) \leq \varepsilon.
$$
For each $\vtau=(\tau(n))_{n \leq N(k)}$, we consider the set
$$
X_{\vtau}=\bigcap_{n \leq N(k)}S_{n,\tau(n)}
$$
so that if $I=\{1,\ldots,N(k)\}$ we have $(X_{\vtau},I,2^{-k}) \in
\D_k$ and thus $\psi(X_{\vtau}) \leq 2^{-k}$, and hence
$\mu(X_{\vtau}) \leq \varepsilon$.

Let us denote by $\Av$ the average over all values of $\vtau$, so
that
\begin{equation}\label{3.9}
\int \Av(1_{X_{\vtau}}(\vz))d\mu(\vz)=\Av \int
1_{X_{\vtau}}(\vz)d\mu(\vz)=\Av \mu(X_{\vtau}) \leq \varepsilon.
\end{equation}
It should be clear that the quantity $\Av(1_{X_{\vtau}}(\vz))$ is
independent of $\vz$. Its value $a_k$ satisfies
$$
a_k=\int \Av 1_{X_{\vtau}}(\vz)d\lambda(\vz)= \Av \int
1_{X_{\vtau}}(\vz) d\lambda(\vz)
$$
where $\lambda$ denotes the uniform measure on $T$. Now
$$
\int 1_{X_{\vtau}}(\vz) d\lambda(\vz)=\lambda(X_{\vtau})=\prod_{n
\leq N(k)}(1-2^{-n})
$$
is bounded below independently of $k$, so that $a_k$ is bounded below independently of $k$.  Finally \eqref{3.9} yields
$$
\varepsilon \geq \int \Av (1_{X_{\vtau}}(\vz))d\mu(\vz)=a_k\mu(T),
$$
and since $\varepsilon$ is arbitrary this shows that $\mu(T)=0$.\qed

\begin{teo}\label{teo3.3} If the sequence $N(k)$ is chosen as in
\eqref{3.7}, for each $\varepsilon >0$ we can find an
$\varepsilon$-exhaustive submeasure $\nu \leq \psi$.
\end{teo}

This result is of course much weaker than Theorem \ref{teo1.1}. We
present its proof for pedagogical reasons. Several of the key
ideas required to prove Theorem \ref{teo1.1} will already be
needed here, and should be much easier to grasp in this simpler
setting.

Given $A \in \A_m$, let us define the map $\pi_A:T \to A$ as
follows: If $\tau_1,\ldots,\tau_m$ are such that
$$
\vz \in A \Longleftrightarrow \forall i \leq m,\ \ z_i=\tau_i
$$
then for $\vz \in T$ we have $\pi_A(\vz)=\vy$ where
$$
\vy=(\tau_1,\ldots,\tau_m,z_{m+1},\ldots).
$$

\begin{df}\label{df3.4}{\bf (Farah)} Given $m<n$, we say that a
set $X\subset T$ is $(m,n,\psi)$-thin if
$$
\forall A \in \A_m,\ \ \exists C \in \B_n, \ \ C \subset A,\ \ C
\cap X=\varnothing,\ \ \psi(\imin {\pi_A}(C)) \geq 1.
$$
\end{df}

The idea is now that in each atom of rank $m$, $X$ has a
$\B_n$-measurable hole that is large with respect to $\psi$. Of
course, we cannot require that $\psi(C) \geq 1$ because $\psi(C)
\leq \psi(A)$ will be small, and one should think of $\psi(\imin
{\pi_A}(C))$ as measuring the ``size of $C$ with respect to
$A$''.

Obviously, if $n' \geq n$ and if $X$ is $(m,n,\psi)$-thin, it is
also $(m,n',\psi)$-thin. For a subset $I$ of $\N^*$, we say that $X$
is $(I,\psi)$-thin if it is $(m,n,\psi)$-thin whenever $m,n \in
I$, $m<n$. By the previous observation, it suffices that this
should be the case when $m$ and $n$ are consecutive elements of
$I$.

Consider a given integer $q$ and consider an integer $b$, to be
determined later. Consider the class $\F$ of marked weighted sets
defined as
$$
\F=\{(X,I,w);\, X \text{ is $(I,\psi)$-thin, $\card I =b,$
$w=2^{-q}$}\}.
$$
We define
$$
\nu=\varphi_{\F \cup \D},
$$
where $\D$ is the class \eqref{2.6}. Thus $\nu \leq \psi=\varphi_\D$, so it
is pathological.
\begin{prop}\label{prop3.5} The submeasure $\nu$ is
$2^{-q}$-exhaustive.
\end{prop}
\begin{prop}\label{prop3.6} If we assume
\begin{equation}\label{3.11}
b=2^{2q+10}
\end{equation}
we have $\nu(T) \geq 2^4.$
\end{prop}

Both these results assume that \eqref{3.7} holds. This condition
is assumed without further mention in the rest of the paper. \\
We first prove Proposition \ref{prop3.5}. Again, the arguments are
due to I. Farah \cite{F} and are of essential importance.

\begin{lem}\label{lem3.6} Consider a  sequence $(E_i)_{i
\geq 1}$ of $\B$ and assume that
$$
\forall n,\ \ \psi\left(\bigcup_{i \leq n} E_i\right) <1.
$$
Assume that for a certain $m \geq 1$, the sets $E_i$ do not depend
on the coordinates of rank $\leq m$. Then for each $\alpha >0$ we
can find a set $C \in \B$, that does not depend on the coordinates
of rank $\leq m$, and satisfies $\psi(C) \leq 2$ and
$$
\forall i \geq 1,\ \ \psi(E_i \setminus C) \leq \alpha.
$$
\end{lem}
\dim\ \ By definition of $\psi$ for each $n$ we can find a finite
set $F_n \subset \D$ with $w(F_n)<1$ and $\bigcup_{i \leq n}E_i
\subset \cup F_n$. For an integer $r \geq m+2$, let
\begin{equation}\begin{split}
F_n^r=\{(X,I,w) \in F_n;\, \card I \cap \{m+1,\ldots,r-1\} < \card
I/2; \\
\card I \cap \{m+1,\ldots,r\} \geq \card I/2\},
\end{split}\end{equation}
so that the sets $F_k^r$ are disjoint as $r$ varies. We use Lemma
\ref{lem3.1} and \eqref{3.5a} with $J=I \cap \{m+1,\ldots,r\}$ to
obtain for each $(X,I,w)$ an element $(X',I',w')$ of $\D$ such
that $X' \supset X$, $w' \leq 2w$, and $X'$ depends only on the
coordinates of rank in $\{m+1,\ldots,r\}$ (or, equivalently, $I \subset \{m+1, \cdots, r\}$). We denote by $F'_n{}^r$ the
collection of the sets $(X',I',w')$ as $(X,I,w) \in F_n^r$. Thus
$\cup F_n'{}^r \supset \cup F_n^r$, and $w(F'_n{}^r) \leq
2w(F_n^r)$.

Consider an integer $i$, and $j$ such that $E_i \in \B_j$. We
prove that for $n \geq i$ we have $E_i \subset \bigcup_{r \leq
j}\cup F_n'{}^r$. Otherwise, since both these sets depend only on
the coordinates of rank in $\{m+1,\ldots,j\}$, we can find a set
$A$ depending only on those coordinates with $A \subset E_i
\setminus \bigcup_{r \leq j}\cup F_n'{}^r$, and thus $A \subset E_i
\setminus \bigcup_{r \leq j} \cup F_n^r$. Since $E_i \subset \cup
F_n$, we have $A \subset \cup F^\sim$, where $F^\sim=F_n \setminus
\bigcup_{r \leq j}  F_n^r$. Now, by definition of $F^r_n$, if
$(X,I,w) \in F^\sim$, we have $\card (I \setminus
\{m+1,\ldots,j\}) \geq \card I /2$. Again use  Lemma
\ref{lem3.1},  now with $J=\{m+1,\ldots,j\}^c$ to see that we can
find $(X',I',w')$ in $\D$ with $w' \leq 2w$ and $X' \supset X$,
$X'$ does not depend on the coordinates of rank in
$\{m+1,\ldots,j\}$. Let $F'$ be the collection of these triplets
$(X',I',w')$, so $F' \subset \D$ and $w(F') \leq 2w(F_n) \leq 2$.
Now $\cup F' \supset \cup F^\sim \supset A$, and since $\cup F'$ does
not depend on the coordinates in $\{m+1,\ldots,r\}$, while $A$ is
determined by these coordinates, we have $\cup F'=T$. But this
would imply that $\psi(T) \leq 2$, while we have proved that
$\psi(T) \geq 2^5$.

Thus $E_i \subset \bigcup_{r \leq j} \cup F_n'{}^r$. For $(X,I,w)$
in $F_n'{}^r$, we have $I \subset \{m+1,\ldots,r\}$. Under
\eqref{3.7} we have  that if $(X,I,w) \in \D_k \cap  F_n'{}^r$ then 
\begin{equation}
\label{e1}
w=2^{-k} \left({N(k) \over \card I}\right)^{\alpha(k)} \geq {2^5
\over {\card I}^{\alpha(k)}} \geq {2^5 \over {r^{\alpha(k)}}},
\end{equation}
 which shows (since $w(F_n^r) \leq 1$)
that $k$ remains bounded independently of $n$. Since moreover $I\subset \{m+1, \cdots, r\}$ there exists a finite set $\D^r \subset \D$ such that $F^r_n \subset \D^r$ for all $n$. Then, by taking a subsequence if necessary, we can assume that for each $r$ the sets $F^r_n$ are eventually equal to a set $F^r$.  For each triplet $(X,I,w)$ in  $F^r$, the set $X$ depends only on the coordinates of
rank in $\{ m+1,\ldots, r\}$, and it should be obvious that $\sum_{r \geq m}w(F^r) \leq 2$ and $E_i
\subset \bigcup_{r \leq j} \cup F^r$ (whenever $j$ is such that
$E_i \in \B_j$).

Consider $r_0$ such that $\sum_{r >r_0} w(F^{r}) \leq \alpha$,
and let $C=\bigcup_{r \leq r_0} \cup F^r$. Thus $C \in \B$, $C$
does not depend on the coordinates of rank $\leq m$ and $\psi(C)
\leq \sum_{r \leq r_0} w(F^r) \leq 2$. Moreover, since $E_i
\subset \bigcup_{r \leq j} \cup F^r$ whenever $j$ is large enough that
$E_i \in \B_j$, we have
$$
E_i \setminus C \subset \bigcup_{r_0<r\leq j} \cup F^r,
$$
so that $\psi(E_i \setminus C) \leq \sum_{r>r_0} w(F^r) \leq
\alpha$. \qed

\begin{lem}\label{lem3.7}{\bf (Farah)} Consider $\alpha >0$, $B
\in \B_m$, and a disjoint sequence $(E_i)$ of $\B$. Then there
exists $n>m$, a set $B' \subset B$, $B' \in\B_n$, so that $B'$ is
$(m,n,\psi)$-thin and $\limsup_{i \to \infty} \psi((B \cap E_i)
\setminus B') \leq \alpha$.
\end{lem}
\dim\ \ Consider $\alpha '=\alpha /\card \A_m$. Consider $A \in
\A_m$, $A \subset B$.\\
\emph{Case 1.} $\exists p$; $\psi\left(\imin
{\pi_A}\left(\bigcup_{i \leq p} E_i \right)\right) \geq 1$.

We set $C'=C'(A)= A \setminus \bigcup_{i \leq p} E_i$, so that $\psi(\imin
{\pi_A}(A \setminus C')) \geq 1$ and $(A \cap E_i) \setminus C'=\varnothing$ for
$i>p$.\\
\emph{Case 2.} $\forall p$; $\psi\left(\imin
{\pi_A}\left(\bigcup_{i \leq p} E_i \right)\right) < 1$.

The sets $\imin {\pi_A}(E_i)$ do not depend on the coordinates of
rank $\leq m$ so by Lemma \ref{lem3.6} we can find a set $C \in \B$, that
does not depend on the coordinates of rank $\leq m$, with $\psi(C)
\leq 2$ and $\limsup_{i \to \infty} \psi\left(\imin {\pi_A}(E_i)
\setminus C\right) \leq \alpha '$. Let $C'=C'(A)=\pi_A(C)\subset A$. Since
$C$ does not depend  on the coordinates of rank $\leq m$, we have
$C=\imin {\pi_A}(C')$ so that $\psi\left(\imin {\pi_A}(C')\right) \leq
2$. Since $\pi_A(\vz)=\vz$ for $\vz \in A$, we have
$$
(A \cap E_i) \setminus C' \subset \imin {\pi_A}(E_i) \setminus C
$$
so that
$$
\limsup_{i \to \infty} \psi((A \cap E_i) \setminus C') \leq
\limsup_{ i \to \infty} \psi \left(\imin {\pi_A}(E_i) \setminus
C\right) \leq \alpha '.
$$
Let us now define
$$B'=\bigcup\{C'=C'(A);\, A \in \A_m,\, A \subset B\},$$
so that
\begin{equation}\label{3.13}
\limsup_{i \to \infty} \psi((B \cap E_i) \setminus B') \leq \sum
\limsup_{i \to \infty} \psi((A \cap E_i) \setminus C') \leq \alpha
' \card \A_m \leq \alpha,
\end{equation}
where the summation is over $A \subset B$, $A \in \A_m$.

Consider $n$ such that $B' \in \B_n$. To prove that $B'$ is
$(m,n,\psi)$-thin it suffices to prove that $\psi\left(\imin
{\pi_A}(A \setminus C')\right) \geq 1$ whenever $A \in \A_m$, $A
\subset B$, because $B' \cap A=C'$, and thus $A \setminus B'=A \setminus C'$. This was already done in case 1. In case 2, we observe
that
$$
\psi \left(\imin {\pi_A}(A \setminus C')\right)=\psi\left(\imin
{\pi_A}(C')^c\right)
$$
and that
$$
2^5 \leq \psi(T) \leq \psi\left(\imin
{\pi_A}(C')\right)+\varphi\left(\imin {\pi_A}(C')^c\right) \leq
2+\psi\left(\imin {\pi_A}(C')^c\right).
$$
\qed

\noindent {\bf Proof of Proposition \ref{prop3.5} (Farah).}
Consider a disjoint sequence $(E_i)_{i \geq 1}$ of $\B$. Consider $\alpha >0$. Starting with $B_0=T$, we use Lemma
\ref{lem3.7} to recursively construct sets $B_\ell \in \B$ and
integers $(n_1,n_2,\ldots )$ such that $B_\ell$ is
$(I_\ell,\psi)$-thin for $I_\ell=\{1,n_1,n_2,\ldots,n_\ell\}$ and
$B_\ell \subset B_{\ell -1}$,
\begin{equation}\label{3.14}
\limsup_{i \to \infty} \psi((E_i \cap B_{\ell -1})\setminus
B_{\ell}) \leq \alpha .
\end{equation}
We have, since $B_0=T$
$$
E_i \setminus B_\ell \subset \bigcup_{m \leq \ell}((E_i \cap
B_{m-1}) \setminus B_m),
$$
and the subadditivity of $\psi$ then implies that
$$
\psi(E_i \setminus B_\ell) \leq \sum_{m \leq \ell} \psi((E_i \cap
B_{m-1}) \setminus B_m)
$$
and thus
\begin{equation} \label{3.15}
\limsup_{i \to \infty} \psi(E_i \setminus B_\ell) \leq \alpha
\ell.
\end{equation}
For $\ell =b$ (or even $\ell =b-1$) (where $b$ is given by \eqref
{3.11}) the definition of $\F$ shows
that $(B_\ell,I_\ell,2^{-q}) \in \F$, and thus $\nu(B_\ell) \leq
2^{-q}$. Since $\nu \leq \psi$, we have
$$
\nu(E_i) \leq \nu(B_\ell)+\psi(E_i \setminus B_\ell) \leq
2^{-q}+\psi(E_i \setminus B_\ell),
$$
and \eqref{3.15} shows that
$$
\limsup_{i \to \infty} \nu(E_i) \leq 2^{-q}+\alpha\ell.
$$
Since $\alpha$ is arbitrary,  the proof is complete.
\qed

We turn to the proof of Proposition \ref{prop3.6}. Considering
$F_1 \subset \F$ and $F_2 \subset \D$, we want to show that
$$
w(F_1)+w(F_2) < 2^4 \Longrightarrow T \not \subset (\cup F_1) \cup
(\cup F_2).
$$
Since $w \geq 2^{-q}$ for $(X,I,w) \in \F$, we have $w(F_1) \geq
2^{-q}\card F_1$, so that $\card F_1 \leq 2^{q+4}$. We appeal to
Lemma \ref{lem2.3} with $s=\card F_1$ and $t=b2^{-q-4}$ (which is
an integer by \eqref{3.11}) to see that we can enumerate
$F_1=(X_\ell,I_\ell,w_\ell)_{\ell \leq s}$ and find sets $J_1
\prec J_2 \prec \cdots \prec J_s$ with $\card J_\ell=t$ and
$J_\ell \subset I_\ell$.

Let us enumerate
\begin{equation} \label{3.16}
J_\ell=\{i_{1,\ell},\ldots,i_{t,\ell}\}.
\end{equation}
An essential idea is that each of the pairs
$\{i_{u,\ell},i_{u+1,\ell}\}$ for $1 \leq u \leq t-1$ gives us a
chance to avoid $X_\ell$. We are going for each $\ell$ to choose
one of these chances using a counting argument. For
\begin{equation}
\label{b1}
\vu=(u(\ell))_{\ell \leq s} \in \{1,\ldots,t-1\}^s,
\end{equation}
we define the set
$$
W(\vu)=\bigcup_{\ell \leq s}
]i_{u(\ell),\ell},i_{u(\ell)+1,\ell}],
$$
where for integers $m<n$ we define $]m,n]=\{m+1,\ldots,n\}$.

We consider the quantity
$$
S(\vu)=\sum \{w;\, (X,I,w) \in F_2,\, \card (I \cap W(\vu)) \geq
\card I/2\}.
$$
We will choose $\vu$ so that $S(\vu)$ is small. Let us denote by
$\Av$ the average over all possible choices of $\vu$. Then, for
any set $I$, by linearity of $\Av$, we have
\begin{eqnarray*}
\Av(\card (I \cap W(\vu))) & = & \sum_{\ell \leq s} \Av(\card (I
\cap ]i_{u(\ell),\ell}, i_{u(\ell)+1,\ell}])) \\
& = & \sum_{\ell \leq s} {1 \over t-1} \card(I \cap
]i_{1,\ell},i_{t,\ell}]) \leq {1 \over t-1}\card I.
\end{eqnarray*}
Thus, by Markov's inequality, we have
$$
\Av(1_{\{\card(I \cap W(\vu)) \geq \card I/2\}}) \leq  {2 \over
t-1}
$$
and, using linearity of average, we get
$$
\Av(S(\vu)) \leq {2 \over t-1} w(F_2) \leq {2^5 \over t-1} \leq
{2^{q+10} \over b}.
$$
Thus, we can find $\vu$ such that $S(\vu) \leq 2^{q+10}/b$. We fix
this value of $\vu$ once and for all. To lighten notation we set
\begin{equation}
\label{3.x}
W=W(\vu);\ \ m_\ell=i_{u(\ell),\ell},\ \ n_{\ell}
=i_{u(\ell)+1,\ell},\ \ W_\ell=]m_\ell,n_\ell]
\end{equation}
so that $W=\bigcup_{\ell \leq s}W_\ell$, and $n_\ell \leq
m_{\ell+1}$ since $n_\ell \in J_\ell$, $m_{\ell+1} \in
J_{\ell+1}$, $J_\ell \prec J_{\ell+1}$.

Let us define
\begin{eqnarray}
F_3 & = & \{(X,I,w) \in F_2;\, \card(I \cap W) \geq \card I/2\}
\label{3.16a} \\
F_4 & = & \{(X,I,w) \in F_2;\, \card (I\cap W) < \card I/2\},
\label{3.16b}
\end{eqnarray}
so that $F_2=F_3 \cup F_4$, and the condition $S(\vu) \leq
2^{q+10}/b$ means that
$$
w(F_3) \leq {2^{q+10} \over b}.
$$
In particular if $(X,I,w) \in F_3$ we have $w \leq 2^{q+10}/b$.
Since $w \geq 2^{-k}$ for $(X,I,w) \in \D_k$ we see that under
\eqref{3.11} we have
\begin{equation}\label{3.17}
(X,I,w) \in \D_k \cap F_3 \Longrightarrow k \geq q.
\end{equation}

Since $s=\card F_1 \leq 2^{q+4}$ and $W=\bigcup_{\ell \leq
s}W_\ell$, if $\card (I \cap W) \geq \card I/2$, there must exist
$\ell \leq s$ with $\card (I \cap W_\ell) \geq 2^{-q-5}\card I$.
This shows that if we define
\begin{equation}\label{3.17a}
F_3^\ell=\{(X,I,w) \in F_3;\, \card (I \cap W_\ell) \geq 2^{-q-5}
\card I\},
\end{equation}
then we have $F_3=\bigcup_{\ell \leq s}F_3^\ell$.

We appeal to Lemma \ref{lem3.1} with $J=W_\ell$, using the fact
that if $k \geq q$ we have
$$(2^{q+5})^{\alpha(k)} \leq 2$$
(with huge room to spare!), to find for each
$(X,I,w) \in F_3^\ell$ a triplet $(X',I',w') \in \D$ with $X
\subset X'$, $w' \leq 2w$, such that $X'$ depends only on the
coordinates of rank in $W_\ell$. Let $F_3'{}^\ell$ be the
collection of these triplets, so that under \eqref{3.11} we have
$$
w(F'_3{}^\ell) \leq 2w(F_3^\ell) \leq 2w(F_3) \leq {2^5 \over b}
\leq{1\over 2}.
$$
We use again Lemma \ref{lem3.1}, this time for $J$ the complement
of $W$, so that $\card (I \cap J) \geq \card I/2$ for $(X,I,w) \in
F_4$, and we can find $(X',I',w) \in \D$ with $w' \leq 2w$, $X'$
contains $ X$ and depends only on coordinates whose rank is not in
$W$. Let $F_4'$ be the collection of these triplets, so that
$w(F_4') \leq 2w(F_4)<2^5$.

Since $\psi(T) \geq 2^5$, we have $T \not\subset \cup F_4'$, so
that we can find $\vz \in T \setminus \cup F_4'$. Since $\cup
F_4'$ depends only on the coordinates whose rank is not in $W$, if
$\vz ' \in T$ is such that $z_i=z'_i$ for $i \notin W$, then $\vz
' \notin \cup F_4'$. To conclude the proof, we are going to
construct such a $\vz '$ that does not belong to any of the sets
$X_\ell$ or $\cup F'_3{}^\ell$. (Thus $\vz '$ will not belong to
$(\cup F_1) \cup (\cup F_2)$.) First, let $A_1 \in \A_{m_1}$ such
that $\vz \in A_1$. Since $X_1$ is $(m_1,n_1,\psi)$-thin, there
exists $C \in \B_{n_1}$, $C \cap X_1= \varnothing$,
$\psi\left(\imin {\pi_{A_1}}(C) \right) \geq 1$. Since
$w(F_3'{}^1) \leq 1/2$, we therefore have $\imin {\pi_{A_1}}(C)
\setminus C' \neq \varnothing$, where $C'=\cup F_3'{}^1$. Since
$C'$ does not depend on the coordinates of rank $\leq m_1$ we have
$C'=\imin {\pi_{A_1}}(C')$, so that $\imin {\pi_{A_1}}(C) \setminus
\imin {\pi_{A_1}}(C') \neq \varnothing$, and hence $C \setminus C'
\neq \varnothing$. Since $C'$ depends only on the coordinates of
rank in $W_1$, we have $C' \in \B_{n_1}$, and since $C \in
\B_{n_1}$, we can find $A' \in \A_{n_1}$ with $A' \subset C
\setminus C'$, so that $A' \cap X_1=\varnothing$ and $A' \cap \cup
F_3'{}^1=\varnothing$. Next, we find $A_2 \in \A_{m_2}$ with $A_2
\subset A'$ such that if $\vy \in A_2$ then
$$
\forall i,\ \ n_1 < i \leq m_2 \Longrightarrow y_i = z_i,
$$
and we continue the construction in this manner.
\qed

\section{The construction}\label{4}

Given an integer $p$, we will make a construction ``with $p$
levels'', and we will then take a kind of limit as $p \to \infty$.
We consider the sequence $\alpha(k)$ as in \eqref{3.3}, and a
sequence $M(k)$ to be specified later. The only requirement is
that this sequence increases fast enough. We recall the class $\D$
constructed in the previous section.

We construct classes $(\E_{k,p})_{k \leq p}$, $(\C_{k,p})_{k \leq
p}$ of marked weighted sets, and submeasures $(\varphi_{k,p})_{k
\leq p}$ as follows. First, we set
\begin{eqnarray*}
\C_{p,p}=\E_{p,p}=\D && \\
\varphi_{p,p}=\varphi_\D=\psi.
\end{eqnarray*}

Having defined $\varphi_{k+1,p}$, $\E_{k+1,p}$, $\C_{k+1,p}$, we
then set
\begin{eqnarray*}
\E_{k,p}& = & \Bigg\{(X,I,w);\ \ X \in \B,\ \ \text{$X$ is
$(I,\varphi_{k+1,p})$-thin},\ \ \card I \leq M(k), \\
&& \ \  w = 2^{-k}\left({M(k) \over \card I}\right)^{\alpha(k)}\Bigg\} \\
\C_{k,p} & = & \C_{k+1,p} \cup \E_{k,p}\\
\varphi_{k,p} & = & \varphi_{\C_{k,p}}.
\end{eqnarray*}

To take limits, we fix an ultrafilter $\U$ on $\N^*$ and we define the class $ \E_k$ of marked weighted sets by
\begin{eqnarray}
\label{d1}
(X,I,w) \in \E_k & \Longleftrightarrow & \{p;\ \ (X,I,w) \in
\E_{k,p}\} \in \U
\end{eqnarray}

Of course, one can also work with subsequences if one so wishes.
It seems plausible that with further effort one might prove that
$(X,I,w) \in \E_k$ if and only if $(X,I,w) \in \E_{k,p}$ for all
$p$ large enough, but this fact, if true, is not really relevant
for our main purpose.

We define
$$
\C_k=\D \cup \bigcup_{\ell \geq k} \E_\ell=\C_{k+1} \cup \E_k;\ \
\nu_k=\varphi_{\C_k};\ \ \nu=\nu_1.
$$

 Let us assume that 
\begin{equation}
\label{c2}
 M(k) \geq 2^{(k+5)/\alpha(k)}.
\end{equation}
 Then if $w < 2^5$ and $ (X,I,w) \in \E_{r,p}$, since
\begin{equation}
\label{d2}
 w = 2^{-r}\left({M(r) \over \card I}\right)^{\alpha(r)}\geq \frac{2^5}{{\card I}^{\alpha(r)}},
\end{equation}
 $r$ remains bounded independently of $p$. It then  follows from \eqref{d1} that if $w <2^5$ we have 
\begin{equation}
\label{c1}
 (X,I,w) \in \C_k \Longleftrightarrow \{ p; (X,I,w) \in \C_{k,p} \} \in \U.
\end{equation}
\begin{teo}\label{teo4.1} We have $\nu(T) >0$, $\nu$ is
exhaustive, $\nu$ is pathological, and $\nu$ is not uniformly
exhaustive.
\end{teo}

 The hard work will of course be to show that $\nu(T) >0$ and that $\nu$ is exhaustive, but the other two claims are easy. Since $\nu \leq \psi$, it follows from Proposition \ref{prop3.2}
that $\nu$ is pathological. It then follows from the Kalton-Roberts theorem that $\nu$ is not uniformly exhaustive. This can also be seen directly by showing  that $\liminf_{n \to \infty}\inf_{\tau \leq 2^n}
\nu(S^c_{n,\tau})>0$. To see this,
consider $I \subset \N^*$, and for $n \in I$ let $\tau(n) \leq 2^n$.
Then
$$
T \subset \bigcup_{n \in I} S^c_{n,\tau(n)} \cup
\left(\bigcap_{n\in I} S_{n,\tau(n)}\right)
$$
so that by subadditivity we have
\begin{eqnarray*}
1 \leq \nu(T) &\leq& \sum_{n \in I}
\nu(S^c_{n,\tau(n)})+\nu\left(\bigcap_{n \in I}
S_{n,\tau(n)}\right)\\
&\leq& \sum_{n \in I}
\nu(S^c_{n,\tau(n)})+\psi\left(\bigcap_{n \in I}
S_{n,\tau(n)}\right).
\end{eqnarray*}

The definition of $\D$ shows that if $\card I=N(1)$, the last term
is $\leq 1/2$, and thus $\sum_{n \in I} \nu(S^c_{n,\tau(n)}) \geq
1/2$. This proves that $\nu$ is not uniformly exhaustive.

\medskip
  It could be of interest to observe that the submeasure $\nu$ has nice invariant properties. For each $n$ it is invariant under any permutation of the elements of $T_n$. It was observed by Roberts \cite{R} that if there exists an exhaustive  submeasure that is not uniformly exhaustive, this submeasure can be found with the above invariance property. This observation was very helpful to the author. It pointed to what should be a somewhat canonical example.

\section{The main estimate} \label{5}

Before we can say anything at all about $\nu$, we must of course
control the submeasures $\varphi_{k,p}$. Let us define
$$
c_1=2^4;\ \ c_{k+1}=c_k2^{2\alpha(k)}
$$
so that since $\sum_{k \geq 1} \alpha(k) \leq 1/2$ we have
\begin{equation}\label{5.1}
c_k \leq 2^5.
\end{equation}

\begin{teo}\label{teo5.1} Assume that the sequence $M(k)$
satisfies
\begin{equation}\label{5.2}
M(k) \geq 2^{2k+10}2^{(k+5)/\alpha(k)}(2^3+N(k-1)).
\end{equation}
Then
\begin{equation}\label{5.3}
\forall p,\, \forall k \leq p,\ \ \varphi_{k,p}(T) \geq c_k.
\end{equation}
\end{teo}

Of course \eqref{5.2} implies \eqref{c2}. It is the only requirement we need on the sequence $(M(k))$.

The proof of Theorem \ref{teo5.1} resembles that of Proposition
\ref{prop3.6}. The key fact is that the class $\E_{k,p}$ has to a
certain extent the property of $\D_k$ stressed in Lemma
\ref{lem3.1}, at least when the set $J$ is not too
complicated.

The following lemma expresses such a property when $J$ is an interval. We recall the
notation $(X)_n$ of \eqref{2.5b}.

\begin{lem}\label{lem5.2} Consider $(X,I,w) \in \E_{k,p}$, $k<p$, and
$m_0 < n_0$. Let $I'=I \cap ]m_0,n_0]$ and $A \in \A_{m_0}$.
Then if $X'=\left(\imin {\pi_A}(X)\right)_{n_0}$ we have
$(X',I',w') \in \E_{k,p}$ where $w'=w(\card I/\card
I')^{\alpha(k)}$.
\end{lem}
\dim\ \ It suffices to prove that $X'$ is $(I',
\varphi_{k+1,p})$-thin. Consider $m,n \in I'$, $m < n$, so that $m_0 < m < n
\leq n_0$. Consider $A_1 \in \A_m$, and set $A_2 =\pi_A(A_1)
\subset A$, so that $A_2 \in \A_m$. Since $X$ is
$(m,n,\varphi_{k+1,p})$-thin, there exists $C \subset A_2$, $C \in
\B_n$, with $C \cap X=\varnothing$, $\varphi_{k+1,p}(\imin
{\pi_{A_2}}(C)) \geq 1$. Let $C'= A_1 \cap \imin
{\pi_{A_2}}(C)$, so that $C' \in \B_n$.

We observe that if a set $B$ does not depend on the coordinates of
rank $\leq m$, we have
$$
\imin {\pi_{A_1}}(B)=B=\imin {\pi_{A_1}}(B \cap A_1).
$$
Using this for $B=\imin {\pi_{A_2}}(C)$, we get that $\imin
{\pi_{A_1}}(C')=\imin {\pi_{A_2}}(C)$, and consequently
$\varphi_{k+1,p}\left(\imin {\pi_{A_1}}(C')\right) \geq 1$.

It remains only to prove that $C' \cap X'=\varnothing$. This is
because on $A_1$ the maps $\pi_A$ and $\pi_{A_2}$ coincide, so
that, since $C' \subset A_1$, we have $\pi_A(C')=\pi_{A_2}(C') \subset
C$ and hence $\pi_A(C') \cap X =\varnothing$. Thus $C'\cap \pi_A^{-1}(X) = \varnothing$ and since $C'\in \B_n$ we have $C'\cap X' =\varnothing$. \qed

Given $p$, the proof of Theorem \ref{teo5.1} will go by decreasing induction over $k$. For
$k=p$, the result is true since by Proposition \ref{prop3.2} we
have $\varphi_{p,p}(T)=\psi(T) \geq 2^5\geq c_k$.

Now we proceed to the induction step from $q+1$ to $q$. Considering
$F \subset \C_{q,p}$, with $w(F) <c_q$, our goal is to show that
$\cup F \neq T$. Since $\C_{q,p}=\C_{q+1,p} \cup \E_{q,p}$ we have
$F=F_1 \cup F_2$, $F_1 \subset \E_{q,p}$, $F_2 \subset
\C_{q+1,p}$. 

Let $F'_2=F_2\cap \bigcup_{k < q} \D_k$. When $(X,I,w) \in \D_k$ we have $w \geq 2^{-k}\geq 2^{-q}$, and thus 
$$
2^{-q} \card F'_2  \leq w(F'_2) \leq w(F) \leq c_q \leq 2^5
$$
so that $\card F'_2 \leq 2^{q+5}$. Also, for $(X,I,w) \in \D_k$ we have $\card I \leq N(k)$, so that if we set 
\begin{equation}
\label{a1}
 I^* = \bigcup\{ I; (X,I,w) \in F'_2\} 
\end{equation}
we have 
\begin{equation}
\label{a2}
\card I^* \leq t' :=2^{q+5}N(q-1).
\end{equation}

When $(X,I,w) \in \E_{q,p}$ we have $w \geq 2^{-q}$.
Thus
$$
2^{-q} \card F_1  \leq w(F_1) \leq w(F) \leq c_q \leq 2^5
$$
and thus $s:=\card F_1 \leq 2^{q+5}$. Also, when $(X,I,w) \in \E_{q,p}$ we have
$$
2^{-q}\left({M(q) \over \card I}\right)^{\alpha(q)} = w \leq
2^5
$$
so that
\begin{equation}\label{5.3a}
\card I \geq M(q) 2^{-(q+5)/\alpha(q)}
\end{equation}
and hence, if
\begin{equation}
\label{e2}
t=2^{q+8}+t'
\end{equation}
under \eqref{5.2} we have $\card I \geq st$ where $s=\card F_1$.
We follow  the proof of Proposition \ref{prop3.6}. We appeal to
Roberts's selection lemma  to enumerate $F_1$ as $(X_\ell, I_\ell, w_\ell)_{\ell \leq s}$ and find sets $J_1
\prec J_2 \prec \cdots \prec J_s$ with $\card J_\ell=t$ and
$J_\ell \subset I_\ell$.
 We then appeal to  the counting argument of Proposition \ref{prop3.6}, but instead of allowing in \eqref{b1} all the values of $u(\ell) \leq t-1$, we now restrict the choice of $u(\ell$) by 
$$
u(\ell) \in U_\ell = \{ u; 1 \leq u \leq t-1, I^* \cap ]i_{u,\ell}, i_{u+1,\ell}]= \varnothing\}. 
$$
We observe that by \eqref{a2} and \eqref{e2} we have $\card U_\ell\geq 2^{q+8}-1$.

 The counting argument then allows us to find $\vu$ such that (since $w(F_2) \leq 2^5$)
$$
S(\vu) \leq {2 \over 2^{q+8}-1} w(F_2)  \leq 2^{-q-1}.
$$
 Using the notation \eqref{3.x} we have thus constructed 
intervals $W_\ell=]m_\ell,n_\ell]$, $\ell \leq s$, with $n_\ell
\leq m_{\ell +1}$, in such a manner that $X_\ell$ is $(m_\ell, n_\ell, \varphi_{q+1,p})$-thin and that  if $F_3$ is defined by
\eqref{3.16a}  we have that 
\begin{equation}\label{5.5}
w(F_3)  \leq 2^{-q-1}\leq {1 \over 4}.
\end{equation}
 Moreover, if $W=\bigcup_{\ell \leq s} ]m_\ell,n_\ell]$ we have ensured that 
$$
(X,I,w) \in F'_2 \Longrightarrow W \cap I = \varnothing,
$$
so that in particular if we define $F_4$ by \eqref{3.16b} we have
\begin{equation}
\label{a3}
 (X,I,w) \in F_4\:, \:(X,I,w) \in \bigcup_{k < q}\D_k \Longrightarrow W \cap I = \varnothing.
\end{equation}

As before, \eqref{5.5} implies that if $(X,I,w) \in \D_k \cap F_3$, then
$k \geq q$. Let us define the classes $F_3^\ell$, $\ell \leq s$ by 
$$
F_3^\ell=\{(X,I,w) \in F_3;\, \card (I \cap W_\ell) \geq 2^{-q-6}
\card I\},
$$
 so that, since $s \leq 2^{q+5}$, we have $F_3=\bigcup_{\ell \leq s} F_3^\ell$.

\begin{lem}\label{lem5.3} Consider $(X,I,w) \in F_3^\ell$ and $A
\in \A_{m_\ell}$. Then we can find $(X',I',w')$ in $\C_{q+1,p}$
with $X' \supset \imin {\pi_{A}}(X)$, $X' \in \B_{n_\ell}$, $w'
\leq 2w$.
\end{lem}
\dim\ \ If $(X,I,w) \in \D$ we have already proved this statement
in the course of the proof of Proposition \ref{prop3.5}, so, since
$\C_{q+1,p}=\D \cup \bigcup_{q+1 \leq r \leq p}\E_{r,p}$, it
suffices to consider the case where $(X,I,w) \in \E_{r,p}$, $r
\geq q+1$. In that case, if $I'=I \cap W_\ell$, we have
$$
\left({\card I \over \card I'}\right)^{\alpha(r)} \leq
(2^{q+6})^{\alpha(r)} \leq 2
$$
and the result follows from Lemma \ref{lem5.2}. \qed

\begin{cor}\label{cor5.4} Consider $A \in \A_{m_\ell}$. Then there
is $A' \in \A_{n_\ell}$ such that $A' \subset A$, $A' \cap
X_\ell=\varnothing$ and $A' \cap \cup F_3^\ell=\varnothing$.
\end{cor}
\dim\ \ Lemma \ref{lem5.3} shows that $\imin {\pi_A}(\cup F_3^\ell) \subset
C' $, where $C' \in \B_{n_\ell}$ and $\varphi_{q+1,p}(C') \leq
2w(F_3^\ell) \leq 1/2$. Since $X_\ell$ is
$(m_\ell,n_\ell,\varphi_{q+1,p})$-thin, we can find $C \in
\B_{n_\ell}$, $C \subset A$, $C \cap X=\varnothing$ with
$\varphi_{q+1,p}\left(\imin {\pi_A}(C)\right) \geq 1$. Thus we
cannot have $\imin {\pi_A}(C) \subset C'$ and hence since both these sets belong to $\B_{n_\ell}$ we can find
$A_1 \in \A_{n_\ell}$ with
$$
A_1 \subset \imin {\pi_A}(C) \setminus C' \subset \imin {\pi_A}(C)
\setminus \imin {\pi_{A}}(\cup F_3^\ell).
$$
Thus $A'=\pi_A(A_1) \in \A_{n_\ell}$,  $A' \cap \cup F_3^\ell=\varnothing$, $A' \subset
C$, so that $A'\cap X_\ell=\varnothing$. \qed

We now construct a map $\Xi:T \to T$ with the following properties. For $\vy \in T$,
$\vz=\Xi(\vy)$ is such that $\vz_i=\vy_i$ whenever $i \notin
W=\bigcup_{\ell \leq s} ]m_\ell,n_\ell]$. Moreover, for each
$\ell$, and each $A \in \A_{m_\ell}$, there exists $A'\in \A_{n_\ell}$ with
$$
\vy \in A \Longrightarrow \Xi(\vy)\in A',
$$
and $A'$ satisfies $A' \cap X_\ell =\varnothing$ and $A' \cap \cup
F_3^\ell=\varnothing$.

The existence of this map is obvious from Corollary \ref{cor5.4}.
It satisfies
\begin{equation}\label{5.6}
\ell \leq s \Longrightarrow \Xi(T) \cap X_\ell=\varnothing,\ \
\Xi(T) \cap \cup F_3^\ell=\varnothing.
\end{equation}

 It has the further property  that
for each integer $j$ the first $j$ coordinates of $\Xi(\vy)$
depend only on the first $j$ coordinates of $\vy$.

We recall that $F_4$ is as in \eqref{3.16b}.

\begin{lem}\label{lem5.5} We have $\varphi_{q+1,p}\left(\imin
\Xi(\cup F_4)\right) < c_{q+1}.$
\end{lem}

\noindent {\bf Proof of Theorem \ref{teo5.1}.} Using the induction hypothesis $\varphi_{q+1,p}(T) \geq c_{q+1}$
we see that there is $\vy$ in $T \setminus \imin \Xi(\cup F_4)$,
so that $\Xi(\vy) \notin \cup F_4$. Combining with \eqref{5.6} we
see that $\Xi(\vy) \notin \bigcup_{\ell \leq s}X_\ell= \cup F_1$, $\Xi(\vy)
\notin \cup F_3$, so that $\Xi(\vy) \notin \cup F$. \qed

\noindent {\bf Proof of Lemma \ref{lem5.5}.} We prove that if
$(X,I,w) \in F_4$, then $\varphi_{q+1,p}\left(\imin \Xi(X)\right)
\leq w2^{2\alpha(q)}$. This suffices since $w(F_4)<c_q$.\\

\emph{Case 1.} $(X,I,w) \in \D_k$, $k<q$.\\

In that case, by \eqref{a3} we have $I \cap W=\varnothing$, so that
$\imin \Xi(X)=X$ and thus $\varphi_{q+1,p}\left(\imin
\Xi(X)\right)=\varphi_{q+1,p}(X) \leq w$.\\

\emph{Case 2.} We have $(X,I,w) \in \D_k$, $k \geq q$.\\

We use Lemma \ref{3.1} with $J=\N^* \setminus W$ and the fact that
$\alpha(k) \leq \alpha(q) \leq (q+5)^{-3}$. This has already been
done in the previous section.\\

\emph{Case 3.} $(X,I,w) \in \E_{r,p}$ for some $q+1 \leq r <
p$.\\

In a first stage we prove the following. Whenever $m,n \in I$ are
such that $m<n$, and $]m,n]\cap W=\varnothing$, then $\imin
\Xi(X)$ is $(m,n,\varphi_{r+1,p})$-thin. Since 
for each integer $j$ the first $j$ coordinates of $\Xi(\vy)$
depend only on the first $j$ coordinates of $\vy$, whenever $A
\in \A_m$ there is $A' \in\A_m$ with $\Xi(A) \subset A'$. Since
$X$ is $(m,n,\varphi_{r+1,p})$-thin we can find $C' \in \B_n$ with
$C' \cap X=\varnothing$, $C' \subset A'$, and
$\varphi_{r+1,p}\left(\imin {\pi_{A'}}(C')\right) \geq 1$. Let
$C=\imin \Xi(C') \in \B_n$. We observe that $C \cap \imin
\Xi(X)=\varnothing$ and we  now prove 
that
\begin{equation}\label{5.8}
\Xi\left(\pi_A\left(\imin {\pi_{A'}}(C')\right)\right) \subset C'.
\end{equation}
 Consider $\tau_1,\ldots,\tau_m$ and
$\tau'_1,\ldots,\tau'_m$ such that
\begin{eqnarray*}
A & = & \{\vz \in T;\ \ \forall i \leq m,\ \ z_i=\tau_i\} \\
A' & = & \{\vz \in T;\ \ \forall i \leq m,\ \ z_i = \tau'_i\}.
\end{eqnarray*}
Consider $\vy \in \imin {\pi_{A'}}(C')$. Then there exists $\vy ' \in C'$
with $y_i=y'_i$ for $i>m$. Thus $\vy''=\pi_A(\vy)$ is such that
$y''_i=\tau_i$ for $i \leq m$, and $y_i''=y'_i$ for $i>m$, so that
$\vz=\Xi(\vy'')$ is such that $z_i=\tau'_i$ for $i<m$. Moreover $z_i=y''_i$ for $i \not \in W$, and since $]m,n]\cap W = \varnothing$, we have 
$z_i=y''_i= y'_i$ for $m<i\leq n$. Since $C' \subset A'$, we have
$y'_i=\tau'_i$ for $i<m$, so that $z_i=y'_i$ for all $i \leq n$,
and thus $\vz \in C'$ since $\vy ' \in C' \in \B_n$. Since $\vy$ is arbitrary this proves
\eqref{5.8}, which implies that
$$
\imin {\pi_{A'}}(C') \subset \imin {\pi_A}\left( \imin
\Xi(C')\right)=\imin {\pi_A}(C),
$$
so that $\varphi_{r+1,p}\left(\imin {\pi_A}(C)\right) \geq 1$ and we
have proved that $\imin \Xi(X)$ is $(m,n,\varphi_{r+1,p})$-thin.

For each $\ell \leq 1$, consider the largest element $i(\ell)$ of
$I$  that is $\leq m_\ell$. (Trivial modifications of the
argument take care of the case where $I$ has no elements $\leq
m_\ell$). Let
$$
I'=I \setminus (W \cup \{i(1),\ldots,i(s)\}),
$$
so that, since $\card (I \setminus W) \geq \card I/2$, we have
$$
\card I' \geq {\card I \over 2} - s \geq {\card I \over 2} -
2^{q+5} \geq {\card I \over 4},
$$
using \eqref{5.3a} and \eqref{5.2}. We claim that $\imin \Xi(X)$ is
$(m,n,\varphi_{r+1,p})$-thin whenever $m<n$, $m,n \in I'$. To see
this, consider the smallest element $n'$ of $I$ such that $m<n'$.
Then $n'\leq n$, so it suffices to show that $\imin \Xi(X)$ is
$(m,n',\varphi_{r+1,p})$-thin. By the first part of the proof, it
suffices to show that $W \cap ]m,n']=\varnothing$. Assuming
$W_\ell \cap ]m,n'] \neq \varnothing$, we see that $m_\ell < n'$.
 Since $m \not \in W_\ell$ we have $m \leq i(\ell)$  and since $m \neq i(\ell)$, we have $m<i(\ell)\leq m_\ell$,
contradicting the choice of $n'$.

Let $w'=w(\card I / \card I')^{\alpha(q)} \leq w2^{2\alpha(q)}$. It should then be
obvious that $(\imin \Xi(X),I',w') \in \E_{r,q}$, so that
$\varphi_{q+1,p}\left(\imin \Xi(X)\right) \leq w2^{2\alpha(q)}$.
\qed

\section{Exhaustivity}\label{6}

\begin{lem}\label{lem6.1} Consider $B \in \B$ and $a>0$. If
$\nu_k(B) < a$ then
$$
\{p;\, \varphi_{k,p}(B)<a\} \in \U.
$$
\end{lem}
\dim\ \ By definition of $\nu_k=\varphi_{\C_k}$, there exists a
finite set $F \subset \C_k= \D\cup \bigcup_{r\geq k} \E_r$ with $w(F)<a$ and $\cup F \supset B$. By
definition of $\E_{r}$, for $(X,I,w) \in \E_r$ we have
$$
\{p;\, (X,I,w) \in \E_{r,p}\} \in \U,
$$
so that since $\C_{k,p}= \D \cup \bigcup_{k \leq r < p} \E_{r,p}$ we have $\{p;\, F \subset \C_{k,p}\} \in \U$ and thus $\varphi_{k,p}(B)
\leq w(F)<a$ for these $p$. \qed

\begin{cor}\label{cor6.2} We have $\nu(T) \geq 16$.
\end{cor}
\dim\ \ By Lemma \ref{lem6.1}, and since $\varphi_{1,p}(T) \geq
c_1=16$ by Theorem \ref{teo5.1}.\\

The next lemma is a kind of converse to Lemma \ref{lem6.1}, that
lies much deeper.

\begin{lem}\label{lem6.3} Consider $B \in \B$ with $\nu_k(B) \geq
4$. Then
$$
\{p;\, \varphi_{k,p}(B) \geq 1\} \in \U.
$$
\end{lem}
\dim\ \ Consider $n$ such that $B \in \B_n$, and assume for
contradiction that
$$
U=\{p;\, \varphi_{k,p}(B)<1\} \in \U.
$$
Thus, for $p \in U$, we can find $F_p \subset \C_{k,p}$ with $B
\subset \cup F_p$ and $w(F_p) \leq 1$. Let
\begin{eqnarray*}
F^1_p & = & \{(X,I,w) \in F_p;\, \card (I \cap \{1,\ldots,n\})
\geq
\card I/2\} \\
F^2_p & = & F_p \setminus F_p^1=\{(X,I,w) \in F_p;\, \card (I
\cap \{1,\ldots,n\}) < \card I/2\}.
\end{eqnarray*}
Using Lemmas \ref{lem3.1} and  \ref{lem5.2} we find a family $F^\sim_p$ of triples
$(X',I',w')$ in $\C_{k,p}$ with $\cup F_p^\sim \supset \cup F_p^1$, $w(F_p^\sim) \leq 2$
and $I' \subset \{1,\ldots,n\}$, $X' \in \B_n$, so that $\cup
F_p^\sim \in \B_n$.

We claim that $B \subset \cup F_p^\sim$. For, otherwise, since $B$
and $\cup F^\sim_p$ both belong to $\B_n$, we can find $ A \in \A_n$ with $A \subset B
\setminus \cup F_p^\sim$, so that $A \subset \cup F_p^2$. By Lemma
\ref{lem5.2} again (or, to be exact, its obvious extension to the case $n_0 = \infty$) and Lemma \ref{lem3.1} we get
$$
\varphi_{k,p}(T)=\varphi_{k,p}\left(\imin {\pi_A}(\cup
F_p^2)\right) \leq 2w(F_p^2) \leq 2,
$$
which is impossible because $\varphi_{k,p}(T) \geq 16$.

 Using \eqref{e1} and \eqref{d2} we see that there exists a finite collection $\cal G$ of triplets $(X,I,w)$ such that $F_p^\sim  \subset \cal G$ for all $p$. Thus there exists a set $F$ such that $\{p \in U; F_p^{\sim} =F\} \in \U$. If follows from \eqref{c1} that $F \subset \C_k$ and it is obvious that $B \subset \cup F$ and $w(F) \leq 2$, so that $\nu_k(B)
\leq 2$, a contradiction.\qed

\begin{cor}\label{cor6.4} Consider a triplet $(X,I,w)$ and $k$
with $\card I \leq M(k)$ and 
$$
w = 2^{-k}\left({M(k) \over \card I}\right)^{\alpha(k)}.
$$
Assume that $X$ is $(I,\nu_{k+1}/4)$-thin, i.e.
\begin{equation}\label{6.2}\begin{split}
\forall m,n \in I,\ \ m<n,\ \ \forall A \in \A_m,\ \ \exists C \in
\B_n,\ \ C \cap X=\varnothing,\ \ \nu_{k+1}\left(\imin
{\pi_A}(C)\right) \geq 4.
\end{split}\end{equation}
Then $(X,I,w) \in \E_k$.
\end{cor}
\dim\ \ If $\nu_{k+1}(\imin {\pi_A}(C)) \geq 4$ then by Lemma \ref{lem6.3} we have $\{p;\,
\varphi_{k+1,p}(\imin {\pi_A}(C)) \geq 1\} \in \U$ and
$$
\{p;\, (X,I,w) \in \E_{k,p}\} \supset \bigcap\left  \{p;\,
\varphi_{k+1,p}\left(\imin {\pi_A}(C)\right) \geq 1\right \} \in \U,
$$
where the intersection is over all sets $A,C$ as in \eqref{6.2}.
\qed

\begin{lem}\label{lem6.5} Consider a sequence $(E_i)$ of
$\B$, and assume that these sets do not depend on the coordinates
of rank $\leq m$ for a certain $m$. Assume that
$$
\forall n,\ \ \nu_k\left(\bigcup_{i \leq n} E_i\right) < 4.
$$
Then for each $\alpha >0$ there is $C \in \B$, that does not
depend on the coordinates of rank $\leq m$, and such that
$\nu_k(C) \leq 8$ and $\nu_k(E_i \setminus C) \leq \alpha$
for each $i$.
\end{lem}
\dim\ \ For each $n$, let
$$U_n=\left\{p;\, \varphi_{k,p}\left(\bigcup_{i \leq n} E_i\right)
< 4\right\}
$$
so that $U_n \in \U$ by Lemma \ref{lem6.1}. For $p \in U_n$ we can
find $F_{n,p} \subset \C_{k,p}$ with $\bigcup_{i \leq n}E_i
\subset \cup F_{n,p}$ and $w(F_{n,p}) \leq 4$. For $r \geq
m+1$ we define
\begin{equation*}\begin{split}
F_{n,p}^r=\left\{ (X,I,w) \in F_{n,p};\, \card (I \cap
\{m+1,\ldots,r-1\}) \leq {1\over 2} \card I;\right. \\
\left.\card (I \cap \{m+1,\ldots,r\}) \geq {1\over 2}\card
I\right\},
\end{split}\end{equation*}
and we define 

$$
F'_{n,p}=\left\{(X,I,w) \in F_{n,p};\, \card (I \cap
\{1,\ldots,m\}) \geq {1 \over 4}\card I\right\}.
$$
We use Lemmas \ref{lem3.1} and \ref{lem5.2} to find a set $B \in \B_m$ with
$\varphi_{k,p}(B) \leq 8$ and $B \supset \cup F'_{n,p}$ so that
since $\varphi_{k,p}(T) \geq 16$ we have $B \neq T$ and thus there
exists $A_{n,p} \in \A_m$ with $A_{n,p} \cap \cup
F'_{n,p}=\varnothing$. We use again Lemmas \ref{lem3.1} and \ref{lem5.2} to see that
for $(X,I,w) \in F_{n,p}^r$ we can find $w' \leq 2w$ such that
$(X',I',w')=((\imin {\pi_{A_{n,p}}}(X))_r,\, I \cap
\{m+1,\ldots,r\},w') \in \C_{k,p}$. We observe that $X'$ does not
depend on the coordinates of rank $\leq m$. Let
$F'{}^r{}_{\hspace{-3mm}n,p}$ be the collection of the sets
$(X',I ',w')$ for $(X,I,w) \in F^r_{n,p}$ so that $w(F'{}^r{}_{\hspace{-3mm}n,p})\leq
2w(F_{n,p}^r)$. We claim that if $E_i \in \B_j$ we have
\begin{equation}\label{6.3}
E_i \subset \bigcup_{r \leq j} \cup F'{}^r{}_{\hspace{-3mm}n,p}.
\end{equation}
Otherwise, since both sets  depend only on the coordinates of
rank $\geq m$ and $\leq j$, and since $A_{n,p} \in \A_m$, we would find $ A \in \A_j$ with $A \subset A_{n,p}$ and $A \subset
E_i \setminus \bigcup_{r \leq j} \cup
F'{}^r{}_{\hspace{-3mm}n,p}$. Since $\imin {\pi_{A_{n,p}}}(X) \cap
A_{n,p} \supset X \cap A_{n,p}$, this shows that $A \subset E_i
\setminus \bigcup_{r \leq j} \cup F_{n,p}^r$. Since $A_{n,p} \cap
\cup F^r_{n,p}=\varnothing$ for $r \leq j$, and since $E_i \subset \cup
F_{n,p}$, we  have $A \subset \cup F''_{n,p}$, where
$$
F''_{n,p}=F_{n,p} \setminus \left(F'_{n,p} \cup \bigcup_{r \leq j}
F_{n,p}^r\right) \subset\left\{(X,I,w) \in F_{n,p};\, \card I \cap
\{j+1,\ldots\}\geq {1\over 4} \card I\right\}.
$$
A new application of Lemmas \ref{lem3.1} and \ref{lem5.2} then shows that $T=\imin
{\pi_A}(A)$ satisfies $\varphi_{k,p}(T) \leq 8$, and this is
impossible. So we have proved \eqref{6.3}.\\

Given $r$, we prove using \eqref{e1} and \eqref{d2} that $ F'{}^r{}_{\hspace{-3mm}n,p} \subset {\cal G}^r$ where ${\cal G}^r$ is finite and does not depend on $n$ or $p$. It should then be clear using \eqref{d1}  how to take
limits as $p \to \infty$, $n \to \infty$ to define for $r \geq m+1$ sets $F^r \subset \C_k$ with 
$ \sum_{r \geq m+1}w(F^r) \leq 8$ such  that $E_i \subset \bigcup_{r
\leq j} \cup F^r$ provided $E_i \in B_j$. The elements of $F^r$ are of the type $(X,I,w)$ where $X$
does not depend on the coordinates of rank $\leq m$, and $X \in
\B_r$.

Consider $r_0$ such that $\sum_{r>r_0} w(F^r) < \alpha$ and let
$C=\bigcup_{r \leq r_0} \cup F^r$. Then $\nu_k(C) \leq \sum_{r
\leq r_0}w(F^r) \leq 8$ and
$$
E_i \setminus C \subset  \bigcup_{r_0 \leq r \leq j} \cup F^r
$$
so that $\nu_k(E_i \setminus C) \leq \alpha$.
\qed

\begin{lem}\label{lem6.6} Consider $k>0$, $\alpha >0$, $B
 \in \B_m$, a disjoint sequence $(E_i)$ of $\B$. Then we can find
 $n>m$, a set $B' \in \B_n$, $B' \subset B$ such that $B' $ is
 $(m,n,\nu_{k+1}/4)$-thin and
 $$\limsup_{i \to \infty}\nu_k((B\cap E_i) \setminus B') \leq \alpha.
 $$
\end{lem}
\dim\ \ Nearly identical to that of Lemma \ref{lem3.7}, using
Lemma \ref{lem6.5}, and since $\nu_k(T) \geq 16$.\\

\noindent {\bf Proof that $\nu$ is exhaustive.} For each $k$ we
show that $\nu$ is $2^{-k}$ exhaustive following the method of
Proposition \ref{prop3.5}, and using that by Corollary
\ref{cor6.4}, if $X$ is $(I,\nu_{k+1}/4)$-thin where $\card
I=M(k)$, then $(X,I,2^{-k}) \in \E_k$, so that $\nu(X) \leq \nu_k(X) \leq 2^{-k}$.
\qed

\section{Proof of Theorems \ref{teo1.2} to \ref{teo1.4}}\label{7}

The simple arguments we present here are essentially copied from
the paper of Roberts \cite{R}, and are provided for the
convenience of the reader.

To prove Theorem \ref{teo1.4}, we simply consider the space $\lo$
of real-valued functions defined on the Cantor set that are
$\B$-measurable, provided with the topology induced by the
distance $d$ such that
\begin{equation}\label{7.1}
d(f,0)=\sup\{\varepsilon;\, \nu(\{|f| \geq \varepsilon\}) \geq
\varepsilon\},
\end{equation}
where $\nu$ is the submeasure of Theorem \ref{teo1.1}. We consider
the $\lo$-valued vector measure $\theta$ given by $\theta(A)=1_A$.
Thus $d(0,\theta(A))=\nu(A)$, which makes it obvious that $\theta$
is exhaustive and does not have a control measure. Let us also
note that $d$ satisfies the nice formula
$$
d(f+g,0)\leq d(f,0)+d(g,0),
$$
as follows from the relation $\{|f+g| \geq
\varepsilon_1+\varepsilon_2\} \subset \{|f|\geq \varepsilon_1\}
\cup \{|g| \geq \varepsilon_2\}$.

We start the proof of Theorem \ref{teo1.2}. We first observe  that
the submeasure $\nu$ of Theorem \ref{teo1.1} is strictly positive,
i.e., $\nu(A)>0$ if $A \neq \varnothing$. This follows from subadditivity and the fact that by construction we have $\nu(A) = \nu(A')$ for $A, A'\in \A_n$ and any $n$. 

 We consider the distance
$d$ on $\B$ given by
$$
d(A,B)=\nu(A \Vartriangle B),
$$
where $\Vartriangle$ denotes the symmetric difference.
It is simple to see that the completion $\widehat{\B}$ of $\B$
with respect to this distance is still a Boolean algebra, the
operations being defined by continuity, and that $\nu$ extends to
$\widehat{\B}$ in a positive submeasure, still denoted by $\nu$.
We claim that $\nu$ is exhaustive. To see this, consider a
disjoint sequence $(E_n)$ in $\widehat{\B}$. Consider
$\varepsilon>0$, and for each $n$ find $A_n$ in $\B$ with $\nu(A_n
\Vartriangle E_n) \leq \varepsilon 2^{-n}$. Let $B_n=A_n\setminus
(A_1 \cup \cdots \cup A_{n-1})$, so that, since $E_n=E_n \setminus
(E_1 \cup \cdots \cup E_{n-1})$ we have
\begin{equation}\label{7.2}
\nu(B_n \Vartriangle E_n) \leq \sum_{m \leq n} \nu(E_m
\Vartriangle A_m) \leq \sum_{m \leq n} \varepsilon 2^{-m} \leq
\varepsilon .
\end{equation}
Since the sequence $(B_n)$ is disjoint in $\B$, we have $\lim_{n
\to \infty}\nu(B_n)=0$, and by \eqref{7.2} we have $\limsup_{n \to
\infty}\nu(E_n) \leq \varepsilon$. As $\varepsilon$ is arbitrary,
this proves the result.

Consider now a decreasing sequence $(A_n)$ of $\widehat{\B}$. The
fundamental observation is that it is a Cauchy sequence for $d$.
Otherwise, we could find $\varepsilon >0$ and numbers $m(k)
<n(k)\leq m(k+1)<n(k+1)\cdots$ with $\nu(A_{n(k)} \setminus
A_{m(k)}) \geq \varepsilon$, and this contradicts exhaustivity.

The limit of the sequence $(A_n)$ in $\widehat{\B}$ is clearly the
infimum of this sequence. This shows that $\widehat{\B}$ is
$\sigma$-complete and that $\nu$ is continuous.

It  follows that $\nu$ is countably subadditive, i.e.
\begin{equation}\label{7.3}
\nu\left(\bigcup_{n \geq 1} A_n\right) \leq \sum_{n \geq 1}
\nu(A_n).
\end{equation}
This is because for each $m$, if we set $B_m=\bigcup_{n \geq 1}A_n
\setminus \left(\bigcup_{1\leq n\leq m}A_n\right)$ we have
$$
\nu\left(\bigcup_{n \geq 1} A_n\right) \leq
\nu\left(B_m \cup \bigcup_{1\leq n \leq m}A_n\right) \leq \nu(B_m)+\sum_{1
\leq n \leq m} \nu(A_n)
$$
and that $\nu(B_m) \to 0$ since $\nu$ is continuous, since the sequence $(B_n)$ decreases and since $\bigcap_{m \geq
1}B_m=0$ (the smallest element of $\widehat{\B}$).

\begin{lem}\label{lem7.1} Consider $A \in \B$, and countable
collections $\C_n$, $n \geq 1$ such that $A \subset \cup \C_n$ for
each $n$. Then for each $\eta>0$ there is $A' \subset A$ with
$\nu(A \setminus A') \leq \eta$ such that for each $n$, $A$ is
covered by a finite subset of $\C_n$.
\end{lem}
\dim\ \ Enumerate $\C_n$ as $(\C_{n,m})_{m\geq 1}$. Since $A
\subset \cup\C_n$, we have $\bigcap_k\left(A \setminus \bigcup_{m
\leq k} C_{n,m}\right)=0$, so that by continuity of $\nu$ there
exists $k(n)$ with $\nu\left(A \setminus \bigcup_{m \leq
k(n)}C_{n,k}\right) \leq \eta 2^{-n}$. The set
$A'=\bigcap_n\bigcup_{m \leq k(n)}C_{n,k}$ is for each $n$ covered
by a finite subset of $\C_n$ and it satisfies $\nu(A \setminus A')
\leq \eta$ by \eqref{7.3}. \qed

Consider  a measure $\mu$ on $\widehat{\B}$. Then $\mu$ is not
absolutely continuous with respect to $\nu$ on $\B$, so that we
can find $\varepsilon >0$ and for each $n$ a set $B_n \in \B$ with
$\nu(B_n) \leq 2^{-n}$ and $\mu(B_n) \geq \varepsilon$. Let $A_n=\bigcup_{m \geq n}B_m$. By
\eqref{7.3} we have $\nu(A_n) \leq \sum_{m\geq n} 2^{-m}\leq
2^{-n+1}$ so that if $A=\bigcap_{n \geq 1}A_n$ we have $\nu(A)=0$
and thus $A=0$. But by monotonicity we have $\mu(A_n) \geq
\varepsilon$, so that $\mu$ is not continuous.

On the other hand, $\nu$ is not absolutely continuous with respect
to $\mu$ on $\B$, so for some $\varepsilon>0$ and each $n$ we can
find $B_n \in \B$ with $\nu(B_n) \geq \varepsilon$ and $\mu(B_n)
\leq 2^{-n}$. Let $A_n =\bigcup_{m \geq n}B_m$ and  $A=\bigcap_{n \geq
1} A_n$, so that $\nu(A_n) \geq \varepsilon$ and $\nu(A) \geq
\varepsilon$ by continuity of $\nu$. We use Lemma \ref{lem7.1} with $\eta=\varepsilon/2$, $\C_n=\{B_m;\,m\geq
n\}$, $A=\bigcap_{n \geq 1}A_n$. We have $\nu(A') \geq \varepsilon
/ 2$ since $\nu(A) \geq \varepsilon$, so that $A' \not = 0$. For each $n$, since
$\mu$ is subadditive, and since $A'$ can be covered by a finite subset of $\C_n$ we have $\mu(A') \leq \sum_{m \geq n}
2^{-m}=2^{-n-1}$. Thus $\mu(A')=0$, and $\mu$ hence is not
positive. This concludes the proof of Theorem \ref{teo1.2}.

To prove Theorem \ref{teo1.3}, we first observe that
$\widehat{\B}$ satisfies the countable chain condition, since
$\nu$ is positive and exhaustive. We prove that it also satisfies
the general distributive law. Given a sequence $(\Pi_n)$ of
partitions and $m \in \N^*$, Lemma \ref{lem7.1} produces a set $C_m$
with $\nu(C^c_m) \leq 2^{-m}$ such that $C_m$ is finitely covered
by every partition $\Pi_n$. And $C_1$, $C_2 \setminus C_1$, $C_3
\setminus (C_1 \cup C_2),\cdots$ is the required
partition. This concludes the proof of Theorem \ref{teo1.3}

\addcontentsline{toc}{section}{

\end{document}